\documentclass[12pt,a4paper]{article}
\usepackage{amsmath}
\usepackage{amsfonts}
\usepackage{amssymb}
\usepackage{makeidx}
\usepackage{graphicx}

\newtheorem{preproof}{{\bf Proof. }\hspace{-.15cm}}

\newcommand{\qs}{\preceq _s}

\newcommand{\s}{\prec_s}

\newcommand{\proofbegin}{\begin{preproof}\rm}
\newcommand{\proofend}{\hfill{$\blacksquare$} \end{preproof}}
\newtheorem{thm}{Theorem}[section]
\newcommand{\theobegin}{\begin{thm}\rm }
\newcommand{\theoend}{\end{thm}}
\newtheorem{lem}[thm]{Lemma}
\newcommand{\lembegin}{\begin{lem}\rm }
\newcommand{\lemend}{\end{lem}}
\newtheorem{prop}[thm]{Proposition}
\newcommand{\propbegin}{\begin{prop}\rm }
\newcommand{\propend}{\end{prop}}
\newtheorem{conj}[thm]{Conjection. }
\newcommand{\conjbegin}{\begin{conj}\rm }
\newcommand{\conjend}{\end{conj}}
\newtheorem{cor}[thm]{Coloroly}
\newcommand{\corbegin}{\begin{cor}\rm }
\newcommand{\corend}{\end{cor}}
\newtheorem{question}[thm]{Question. }
\newcommand{\questionbegin}{\begin{question}\rm }
\newcommand{\questionend}{\end{question}}
\newtheorem{defin}[thm]{definition}
\newcommand{\defbegin}{\begin{defin}\rm }
\newcommand{\defend}{\end{defin}}
\newcommand{\examplebegin}{\begin{preexample}\rm}
\newcommand{\exampleend}{\end{preexample}}

\title{
\baselineskip = 0.8cm
\vskip 1cm
\bf On Maximum Signless Laplacian Estrada Indices of Graphs with Given Parameters}
\author{
\bf H.R.Ellahi$^a$, R.Nasiri$^a$, G.H.Fath-Tabar$^b$, A.Gholami$^a$\\
\small  $^a$ Department of Mathematics, Faculty of Science, University of Qom,\\
\small Qom 37161-46611, I. R. Iran\\
\small $^b$ Department of Mathematics, Faculty of Science, University of Kashan,\\ 
\small Kashan 87317-51167, I. R. Iran
}
\date{}
\begin{document}

\maketitle
\begin{abstract}
\baselineskip=0.5cm
\noindent Signless Laplacian Estrada index of a graph $G$, defined as $SLEE(G)=\sum^{n}_{i=1}e^{q_i}$, where $q_1,q_2,\cdots,q_n$ are the eigenvalues of the matrix $\mathbf{Q}(G)=\mathbf{D}(G)+\mathbf{A}(G)$.
We determine the unique graphs with maximum signless Laplacian Estrada indices among the set of graphs with given number of cut edges,  pendent vertices, (vertex) connectivity and edge connectivity.

 \vskip 3mm



\end{abstract}

\section{Introduction}\label{sec01}

Throughout   this paper, each graph, say $G$, is simple with vertex set $V(G)$ and edge set $E(G)$, such that $|V(G)|=n$.
Let $\mathbf{A}(G)$ and $\mathbf{D}(G)$ denote the adjacency matrix and diagonal matrix of vertex degrees of $G$, respectively.
The (resp. signless)  Laplacian matrix of $G$ denoted by 
$\mathbf{L}(G)=\mathbf{D}(G)-\mathbf{A}(G)$ (resp. $\mathbf{Q}(G)=\mathbf{D}(G)+\mathbf{A}(G)$),see \cite{Merris01, Cvetkovic02}.
We denote the eigenvalues of $\mathbf{A}(G)$, $\mathbf{L}(G)$ and $\mathbf{Q}(G)$ by
 $\lambda_1, \lambda_2,\cdots, \lambda_n$, 
$\mu_1,\mu_2, \cdots, \mu_n$, and $q_1,q_2,\cdots, q_n$, respectively.

First time, Estrada \cite{Estrada01} defined a graph-spectrum-based invariant, named \emph{Estrada index}, as follows:
$$EE(G)=\sum^{n}_{i=1}e^{\lambda_i}$$
the concept of Estrada index has found successful applications in biochemistry and in complex network \cite{Estrada01, Estrada02, Estrada03, Estrada04, Estrada05, Estrada06, Estrada07}.
Furthermore, it has been immensely studied in mathematics: 
In some papers, estimaiting and finding good lower and uper bounds of Estrada index in terms of some parameters of graphs have been discussed \cite{Pena01, Gutman04, Gutman05, Gutman06},
 and in some other papers,  the unique graphs having extremum Estrada index in several subcategories of graphs has been investigated  \cite{Das01,Deng01, Du02,  Zhao01, Zhang01}.
\\
Fath-Tabar et al. \cite{FathTabar01} proposed the \emph{Laplacian Estrada index}, in full analogy with estrada index as
$$LEE(G)=\sum^{n}_{i=1}e^{\mu_i}.$$
Surving lower and uper bounds for LEE in terms of different parameters of graphs, and finding graphs with extremum values of LEE in subcategories of graphs, is a part of reaserches about Laplacian Estrada index. For details see \cite{Bamdad01, FathTabar01,Khosravanirad01, Li01, Zhou01}.

Ayyaswamy et al. \cite{Ayyaswamy01} defined the \emph{signless Laplacian Estrada index} as
$$SLEE(G)=\sum^{n}_{i=1}e^{q_i}.$$
They also established lower and upper bounds for SLEE in terms of the number of vertices and edges.

Note that the Laplacian and signlees Laplacian spectra of bipartite graphs coincide \cite{Grone01,Grone02}. 
Thus, for a bipartite graph $G$, $SLEE(G)=LEE(G)$.
Chemically, since the vast majority of molecular graphs are bipartite, we can use the provided statements in  $SLEE$ for $LEE$, and the interesting case is when $SLEE$ and $LEE$ differ, e.g., fullerences, fluoranthenes and other non-alternant conjugated species \cite{Balaban01, Dias01, Gutman01, Gutman02, Gutman03}.

The paper is organized as follows. In section 3, we provide some lemmas to compaire $SLEE$  of trasformated graphs.
In section 4 and 5, we charachterize the graphs with extremal signless Laplacian Estrada indices, including the unique graphs on $n$ vertices with maximum $SLEE$ among the set of all graphs with given number of cut edges, pendent vertices, (vertex) connectivity and edge connectivity.

\section{Preliminaries}\label{sec02}
Denote by $T_k(G)$ the $k$-th signless Laplacian spectral moment of the graph $G$, i.e., $T_k(G)=\sum^{n}_{i=1}q^{k}_{i}$.
By use of  the Taylor expansions of the function $e^x$, we will given arise to the formula:
\begin{equation}\label{eq001}
SLEE(G)=\sum_{k\geq 0}\frac{T_k(G)}{k!}.
\end{equation}
Moreover, by the following definition and theorem, we can easily compare the SLEE of a graph  and another one.

\defbegin\cite{Cvetkovic01}
A \emph{semi-edge walk}  of length $k$ in a graph $G$ is an alternating sequence $W=v_1 e_1 v_2 e_2 \cdots  v_k e_k v_{k+1}$ of vertices $v_1, v_2, \cdots , v_k, v_{k+1}$ and edges $e_1, e_2, \cdots , e_k$ such that  the vertices $v_i$ and $v_{i+1}$ are end-vertices (not necessarily distinct) of edge $e_i$, for any $i=1, 2, \cdots , k$.
 If $v_1=v_{k+1}$, then we say $W$ is a \emph{closed semi-edge walk}.
\defend
\theobegin\label{th001}\cite{Cvetkovic01}
The signless Laplacian spectral moment $T_k$ is equal to the number of closed semi-edge walks of length $k$.
\theoend
Let $G$ and $G'$ be two graphs, and $x,y\in V(G)$, and $x',y'\in V(G')$.
Denoting by $SW_k(G;x,y)$ the set of all semi-edge walks of length $k$ in $G$, which are starting at vetex $x$, and ending to vertex $y$. 
For convenience, we may denote  $SW_k(G;x,x)$ by $SW_k(G;x)$, and set $SW_k(G)=\bigcup_{x\in V(G)}SW_k (G;x)$.
\\
We use the notation $(G;x,y)\qs(G';x',y')$ for, if $ |SW_k (G;x,y)|\leq|SW_k(G';x',y')|$, for any $k\geq 0$.
Moreover, if  $(G;x,y)\qs(G';x',y')$,
 and there exists some $k_0$ such that $|SW_{k_0}(G;x,y)|<|SW_{k_0} (G';x',y')|$, then we write $(G;x,y)\s(G';x',y')$.

Indeed, by these notations, theorem \ref{th001} will change to the formula:
\begin{equation}\label{eq002}
T_k (G)=|SW_k (G;x)|=|\bigcup_{x\in V(G)}SW_k (G;x)|=\sum_{x\in V(G)}|SW_k (G;x)|
\end{equation}
\section{Lemmas}
The next result immediately follows from eq. \ref{eq001} and eq. \ref{eq002}.
\lembegin\label{p006}
Let $G$ be a graph. If $e$ is an edge such that  $e\not\in E(G)$, Then $SLEE(G)<SLEE(G+e)$.
\lemend
\lembegin\label{p001}
Let $G$ be a graph and $u,v\in V(G)$. If $v$ is a pendent vertex attached to $u$, then  $ (G;u)\qs(G;v)$, with equality if and only if  $deg_G⁡(u)=deg_G⁡(v)=1$.
\lemend
\proofbegin 
The case $k=0$ is trivial. 
Let $k>1$ and $W=veW_1 ev\in SW_k (G;v)$, where $W_1$ is a semi-edge walk of length $k-2\geq 0$ in $G$.
We may construct an injection $f_k:SW_k (G;v)\to SW_k (G;u)$, by $f_k (W)=ueW_1 eu$. Thus  $|SW_k (G;v)|\leq|SW_k (G;u)|$, for any $k\geq 2$.
Moreover, if $deg_G(u)>1$, then  we have $|SW_1 (G;v)|=deg_G⁡(v)=1<deg_G⁡(u)=|SW_1 (G;u)|$.
\\
Note that if $deg_G(u)=1$, then $G$ has an automorphism, interchanging $u$ and $v$. 
\proofend
\lembegin\label{p002}
Let $H_1$ and $H_2$ be two graphs with $u,v\in V(H_1 )$ and $w\in V(H_2)$.
 Let $G_u$ ($G_v$, respectively) be the graph obtained from $H_1$ and $H_2$ by identifying $u$ ($v$, respectively) with $w$. 
If $(H_1;v)\s(H_1;u)$,
 then  $SLEE(G_v)<SLEE(G_u)$  (See figure 1).
\lemend 
\begin{figure}[h]
\center
\includegraphics[scale=0.15]{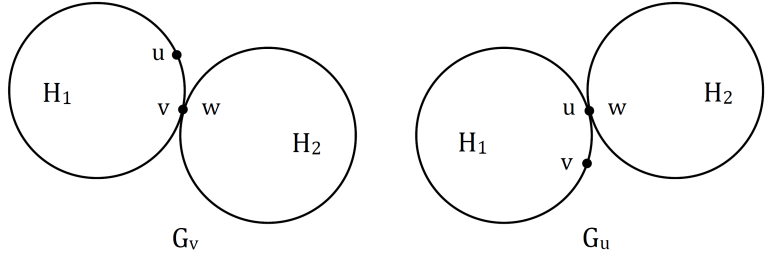}
\caption{A demonstration of the graphs in lemma \ref{p002}}
\end{figure}
\proofbegin
It is enough to show that $T_k (G_v )\leq T_k (G_u)$, for each $k\geq 0$ and there exists a positive integer $k_0$ such that non-equality is strict.

Since $V(G_z)=V(H_1)\cup (V(H_2)\setminus \{w\})$, for each $z\in \{u, v\}$, we have:
$$T_k (G_z )=|\bigcup_{x\in V(H_1)}SW_k (G_z;x)|+|\bigcup_{y\in V(H_2)\setminus \{w\}}SW_k (G_z;y)|.$$

Let $W\in \bigcup_{x\in V(H_1)}SW_k(G_v;x)$.
we can decompose $W$ in a unique form $W=W_1W_2W_3\cdots W_{r-1}W_r$, 
such that $W_1\in  SW_{k_1}(H_1;x,v)$ where $k_1\geq 0$,
and  $W_r\in  SW_{k_r}(H_1;v,x)$ where $k_r\geq 0$,
and when $1<i<r$, we have $k_i>0$, and if  $i$ is even, then $W_i\in  SW_{k_i}(H_2;v)$,
and if $i$ is odd, then $W_i\in  SW_{k_i}(H_1;v)$.

Since  for each $k\geq 0$, $|SW_k (H_1;v)|\leq |SW_k (H_1;u)|$, 
we may consider injections $f_k:SW_k (H_1;v)\to SW_k (H_1;u)$.
Note that $W'=W_rW_1\in SW_{k'}(H_1;v)$, where $k'=k_1+k_r$.
 Thus there is $x'\in H_1$ such that $f_{k'}(W')=W^{'}_{r}W^{'}_{1}$ 
where $W^{'}_{1}\in SW_{k_1}(H_1;x',u)$, and $W^{'}_{r}\in SW_{k_r}(H_1;u,x')$. 

Now, for any $k\geq 0$, we can construct a map:
$$g_k:\bigcup_{x\in V(H_1 )}SW_k (G_v;x)\to \bigcup_{x\in V(H_1 )}SW_k (G_u;x)$$
by
$$g_k(W)=W^{'}_{1}W_2f_{k_3}(W_3)W_4\cdots f_{k_{r-2}}(W_{r-2})W_{r-1}W^{'}_{r}$$
 Indeed, we just replacing any semi-edge walk in $H_1$ by use of injections $f_k$, and fixing another semi-edge walks which are in $H_2$.
The uniqueness of  decomposation of $W$ and $W'$, and being injection of $f_k$ , for any $k\geq 0$, imply that $g_k$ is injective, for any $k\geq 0$.
Therefore
$$|\bigcup_{x\in V(H_1 )}SW_k (G_v;x)|\leq |\bigcup_{x\in V(H_1 )}SW_k (G_u;x)|$$

Similarly, for any $W\in \bigcup_{y\in V(H_2 )\setminus \{w\}}SW_k (G_v;y)$, 
we can decompose $W$ in a unique form $W=W_1W_2W_3\cdots W_{r-1}W_r$, 
where $W_1\in  SW_{k_1}(H_2;y,v)$, where $k_1> 0$ (Note that $y\neq v$),
and  $W_r\in  SW_{k_r}(H_2;v,y)$, where $k_r> 0$,
and when $1<i<r$, we have $k_i>0$, and if  $i$ is odd, then $W_i\in  SW_{k_i}(H_2;v)$,
and if $i$ is even, then $W_i\in  SW_{k_i}(H_1;v)$.
\\
For any $k\geq 0$, we construct a map:
$$h_k:\bigcup_{y\in V(H_2 )\setminus \{w\}}SW_k (G_v;y)\to \bigcup_{y\in V(H_2 )\setminus \{w\}}SW_k (G_u;y)$$
by
$$h_k(W)=W_1f_{k_2}(W_2)W_3\cdots W_{r-2}f_{r-1}(W_{r-1})W_r$$
By the same reasons we said for $g_k$, $h_k$ is an injective map, for any $k\geq 0$, and 
$$|\bigcup_{y\in V(H_2 )\setminus \{w\}}SW_k (G_v;y)|\leq |\bigcup_{y\in V(H_2 )\setminus \{w\}}SW_k (G_u;y)|$$
\\Therefore, $T_k(G_v)\leq T_k(G_u)$, for any $k\geq 0$

To complete the proof, note that for some integer $k_0$ we have $ |SW_{k0} (H_1;v)|<|SW_{k0} (H_1;u)|$.
 It implies that $f_{k0}$ is not surjective map and therefore, there is a closed semi-edge walk $W_0$ in $H_1$ with length $k_0$ started at $u$ which is not covered by $f_{k0}$(and hence, $g_{k_0}$).
Thus, $T_{k_0}(G_v)\leq T_{k_0}(G_u)$. Therefore $SLEE(G_v)< SLEE(G_u)$.
\proofend

\lembegin\label{p003}
Let $G_1$ and $G_2$ be two  graphs with $u\in V(G_1)$ and $v\in V(G_2)$.
 Let $G$ be the graph obtained from $G_1$ and $G_2$, by attaching $u$ to $v$ by an edge,
  and $G'$ be the graph obtained from $G_1$ and $G_2$, by identifying $u$ with $v$, and attaching a pendent vertex to $u$.
  If $deg_G(u),deg_G(v)\geq 2$, then $SLEE(G)<SLEE(G')$.
\lemend
\begin{figure}[h]
\center
\includegraphics[scale=0.15]{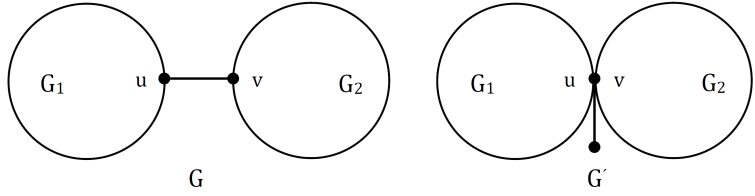}
\caption{A demonstration of graphs in lemma. \ref{p003} ($SLEE(G)<SLEE(G')$).}
\end{figure}
\proofbegin
Let $H_1$ be the graph obtained from $G$ by removing the vertices different from $v$ in $G_2$, and $H_2$ be a copy of $G_2$.
 By applying lemma \ref{p001} we obtained that  $(H_1;v)\s (H_1;u)$.
\\
Now, by applying lemma \ref{p002} on $H_1$ and $H_2$, and assuming $G$ as $G_v$ and $G'$ as $G_u$, we arrive to  $SLEE(G)<SLEE(G')$.
\proofend
\lembegin\label{p009}
Let $G$ be a graph and $v, u, w_1, w_2, \cdots , w_r\in V(G)$.
suppose that
$E_v=\{e_1=vw_1, \cdots , e_r=vw_r\}$ and 
$E_u=\{e'_1=uw_1, \cdots ,  e'_r=uw_r\}$ are subsets of edges, that are not in $G$ (i.e. $e_i,e'_i\not\in E(G)$, for $i=1, 2, \cdots , r$).
Let $G_u=G+E_u$ and $G_v=G+E_v$.
If $(G;v)\s (G;u)$, and $(G;w_i,v)\qs (G;w_i,u)$ for each $i=1,2,\cdots,r$,
Then $SLEE(G_v)<SLEE(G_u)$.
\lemend
\begin{figure}[h]
\center
\includegraphics[scale=0.15]{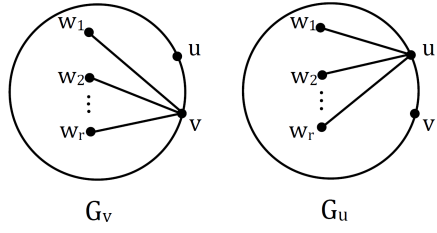}
\caption{An illustration of the graphs  $G_u$ and $G_v$ in lemma \ref{p009}.}
\end{figure}
\proofbegin
Since $(G;v)\prec _s (G;u)$, there exists an injection 
$$f_k:SW_k(G;v)\to SW_k(G,u)$$ for each $k\geq 0$.
Similarly, $(G;w_i,v)\qs (G;w_i,u)$ implies that there exists an injection 
$$f^i_k:SW_k(G;w_i,v)\to SW_k(G,w_i,u)$$
for each $i=1,2,\cdots, r$, and $k\geq 0$.
 While $|SW_k(G,x,y)|=|SW_k(G,y,x)|$ for any $x,y\in V(G)$ (by reversing the semi-edge walk), there exists an injection
$$g^i_k:SW_k(G;v,w_i)\to SW_k(G,u,w_i)$$
for each $i=1,2,\cdots, r$, and $k\geq 0$.

To prove the statement, it is enough to show that $T_k(G_v)\leq T_k(G_u)$,
 and there exists $k_0$ such that inequality is strict.
 Suppose $W\in SW_k(G_v)$.
 we may decompose $W$ to $s+1$ sections $W=W_1e_{j_1}W_2e_{j_2}W_3\cdots W_se_{j_s}W_{s+1}$, where each $W_i$ is a semi-edge walk of length $k_i\geq 0$, in $G$. Obviousely, this decomposition is unique.
 
 Let $1<i<s$. For each $W_i$ one of the following cases hapends:
 
 Case 1. $k=0$ and $W_i=x$, where $x\in \{v,w_1,\cdots, w_r\}$. In this case, we set $h(W_i)=u$, if $x=v$, and $h(W_i)=W_i$, if $x\neq v$.
 
 Case 2. $W_i\in SW_{k_i}(G;v)$. In this case we set $h(W_i)=f_{k_i}(W_i)$.
 
 Case 3. $W_i\in SW_{k_i}(G;w_l,v)$. In this case we set $h(W_i)=f^l_{k_i}(W_i)$.
 
 Case 4. $W_i\in SW_{k_i}(G;v,w_l)$. In this case we set $h(W_i)=g^l_{k_i}(W_i)$.
 
 Case 5. $W_i\in SW_{k_i}(G;w_l,w_j)$, where $l,j=1,2,\cdots,r$. here we set $h(W_i)=W_i$. 
 
 Now, we have one step more, to construct a well-defined injection $h:SW_k(G_v)\to SW_k(G_u)$.
  Since $W$ is closed, $W'=W_{s+1}W_1$ is a semi-edge walk of length $k'=k_{s+1}+k_1$, in $G$, which is in one of the above $5$ cases. 
 Set $W''=h(W')$. We can uniquely decompose $W''=W''_{s+1}W''_1$, where $W''_{s+1}$ is a semi-edge walk of length $k_{s+1}$ in $G$,
 started at $x\in \{u,w_1,\cdots,w_r\}$ and ended at $x'$, and 
  $W''_1$ is a semi-edge walk of length $k_1$ in $G$,
 started at $x'$and ended at  $y\in \{u,w_1,\cdots,w_r\}$.
 
 Finally, it is easy to check that the map 
$h_k:SW_k(G_v)\to SW_k(G_u)$, by 
\begin{align*}
h_k(W)&=h_k(W_1e_{j_1}W_2e_{j_2}W_3\cdots W_se_{j_s}W_{s+1})\\
&=W''_1e'_{j_1}h(W_2)e'_{j_2}h(W_3)\cdots h(W_s)e'_{j_s}W''_{s+1}
\end{align*}
is an injection. Hence, for any $k\geq 0$, $T_k(G_v)\leq T_k(G_u)$.
\\
Moreover, for some $k_0$, $|SW_{k_0}(G,v)|<|SW_{k_0}(G,u)|$, implies that $T_{k_0}(G_v)< T_{k_0}(G_u)$ (Note that $f_{k_0}$ is not surjection). Therefore, $SLEE(G_v)<SLEE(G_u)$.
\proofend
\section{The graph with maximum SLEE with given number of cut edges, and number of pendent vertices}\label{sec03}
Let $a, b \geq 1$. We denote  the set of all graphs which obtained by attaching $b$ pendent vertices to some vertices of $K_a$, by  $\mathcal{G}(a, b)$.
Denote by $G_{a+b,b}$ the graph obtained by attaching $b$ vertices to one vertex of $K_a$,
where $b\geq 0$,  and $K_n$ is the complete graph with $n$ vertices.
\lembegin\label{p004}
Let $a\geq 3$ and $b\geq 1$, and $G\in \mathcal{G}(a, b)$.
 Let $u$ and $v$ be  two distinct non-pendent vertices in $G$. If $u$ has $s$ pendent neighbors in $G$,
and $v$ has $r$ pendent neighbors in $G$,  where $0\leq r<s$,
 Then $(G;v)\s (G;u)$.
\lemend
\proofbegin
Let $V_v=\{v,x_1,x_2,\cdots,x_r\}$, and $E_v=\{e_i=vx_i: 1\leq i\leq r\}$, 
where $x_i$ is a pendant neighbor of $v$, for $i=1,2,\cdots, r$.
 Similarly, let   $V_u=\{u,y_1,y_2,\cdots,y_s\}$, and $E_u=\{e'_i=uy_i: 1\leq i\leq s\}$, 
where $y_i$ is a pendant neighbor of $u$, for $i=1,2,\cdots, s$.

let $W\in SW_k(G;v)$.
$W$ can be decomposed uniquely to $W=W_1W_2W_3$, 
where $W_1$ and $W_3$ are as long as possible,
  consisting of  vertices just in $V_v$, and edges just in $E_v$,
and $W_2$ is begining and ending at some vertices not in $E_v$ 
 (Note that $W_2$ and $W_3$ may be of length $0$).
Set $f_k(W)=W'_1W_2W'_3$, where $W'_i$ is obtained from $W_i$, by replacing vertex $v$ by $u$, and vertices $x_j$ by $y_j$, and edges $e_j$ by $e'_j$, for $i=1,3$, and $j=1,2,\cdots,r$. 
Obviously, $f_k(W)\in SW_k(H;x_1)$, and $f_k:SW_k(H;y_1)\to SW_k(G;u)$ is injective. Thus $|SW_k(G;v)|\leq  |SW_k(G;u)|$,
for $k\geq 1$. Now,    $(G;v)\s (G;u)$ follows from
$|SW_1(G;v)|=deg_G(v)=a-1+r<a-1+s=deg_G(u)=|SW_1(G;u)|. $
\proofend

\lembegin\label{p005}
Let $a\geq 3$, and $b\geq 2$. If $G\in \mathcal{G}(a, b)$, 
then $SLEE(G)\leq SLEE(G_{a+b,b})$,
with equality, if and only if  $G= G_{a+b,b}$.
\lemend
\proofbegin
susppose that $u$ is a vertex of $G$ which has at least one pendent neighbor. 
Since $G\neq G_{a+b,b}$ there is another vertex $v$ of $G$ which has $r\geq 1$ pendent neighbors $w_1,w_2,\cdots ,w_r$.
Let $H_1$ be the graph obtained from $G$ by deleting $w_1,w_2,\cdots ,w_r$,
and $H_2$ be a copy of star $S_{r+1}$,
and  $G'$ be the graph obtained from $H_1$ and $H_2$ by identifying $u$ and the center of star $H_2$. 
Note that $G$ can be obtained from $H_1$ and $H_2$ by identifying $v$ and the center of star $H_2$. 

It follows from lemma \ref{p004} that  
$(H_1;v)\s (H_1;u)$.
Applying lemma \ref{p002} for $G$ as $G_v$, and $G'$ as $G_u$, implies   $SLEE(G)<SLEE(G')$.
By repeating this form of  transformation opration,  and attaching all pendent vertices to $u$, we may finally have $SLEE(G)<SLEE(G_{a+b,b})$.
\proofend
In a connected graph, a \emph{cut edge} is an edge whose removal disconnect the graph. We denote the set of all connected graphs with $n$ vertices  and $r$ cut edges, by $\mathbb{G}(n,r)$, where $0\leq r \leq n-3$.
\theobegin\label{p007}
Let $0\leq r\leq n-3$. If $G\in \mathbb{G}(n,r)$, Then $SLEE(G)\leq SLEE(G_{n,r})$, with equality  if and only if $G= G_{n,r}$.
\theoend

\proofbegin
If $r=0$, then by lemma \ref{p006}, $G_{n,0}=K_n$ has maximum SLEE.

Let $r>1$, and $G$ be a graph in $\mathbb{G}(n,r)$ with maximum SLEE, and $E$ be the set of cut edges in $G$. 
By lemma \ref{p006}, $G-E$ consists of $r+1$ connected components, which are complete. 

If there exists some edge $e$ of $E$, attaching vertices $u$ and $v$ in $G$, 
where $deg_G(u), deg_G(v) \geq 2$,
 then by applying lemma \ref{p003}, we may get a graph in  $\mathbb{G}(n,r)$ with a larger SLEE, a contradiction. 
Therefore, 
there is exactly one end-vertex with degree one for each of edges in $E$, i.e. 
 every cut edge of $G$ has a pendent vertex as an end-vertex. 
Thus $G$ is a graph obtained from $K_{n-r}$,
 by attaching $r$ pendent vertices to some of it's vertices, 
which means  $G\in \mathcal{G}(n-r,r)$.

Now, if $2\leq r\leq n-3$, then by lemma \ref{p005} we have $G= G_{n,r}$. Finally, it is obvious that if $r=1$, $G= G_{n,1}$.
\proofend

Now, we can easily find the unique graph with maximum $SLEE$ among all graphs with $r$ pendent vertices as follows:
\theobegin
Let $0\leq r\leq n-1$.
Among all graphs on $n$ vertices with  $r$ pendent vertices, $G_{n,r}$ is the unique graph which has maximum $SLEE$.
\theoend
\proofbegin\label{p020}
Let $G$ be a graph with $r$ pendent vertices, and have maximum $SLEE$. Let $H$ be the graph obtained from $G$, by removing all pendent vertices.
By lemma \ref{p006}, $H$ is complete graph on $n-r$ vertices. Thus $G\in  \mathcal{G}(n-r,r)$. The cases $r=0,1$ are trivial.  if $n-r\geq 3$, the result follows by lemma \ref{p005}.

If $n-r\leq 2$, then $H=P_{n-r}$. The case $r=n-1$ is trivial. Let $r=n-2$, and $H=P_2$. In this case, the result follows by one time use of lemma \ref{p001}, and lemma \ref{p002}.
\proofend

\section{The Graph with maximum SLEE with given number of vertex  connectivity, and edge connectivity}\label{sec04}

Let $G\cup H$ denote the vertext-disjoint union of graphs $G$ and $H$, and $G\vee H$ be the graph obtained from $G$ and $H$, by attaching any vertex of $G$ to any vertex of $H$. Note that $SLEE(G\cup H)=SLEE(G)+SLEE(H)$.
Let $K_{(p,q)r}=(K_p\cup K_q)\vee K_r$, where $p\geq q\geq 1$, and $r\geq 1$.
\lembegin\label{p012}
If $p\geq q\geq 2$, and 
$r\geq 1$, then 
$SLEE(K_{(p,q)r}) < SLEE(K_{(p+q-1,1)r}).$
\lemend
\proofbegin
Suppose that $V(K_p)=\{x_1,x_2,\cdots, x_p\}$, and $V(K_q)=\{y_1,y_2,\cdots,y_q\}$, and $V(K_r)=\{z_1,z_2,\cdots,z_r\}$.
Let $H$ be the graph obtained from $K_{(p,q)r}$ by removing edges of $y_1$ in $K_q$, i.e. $H=K_{(p,q)r}-\{y_1y_i:2\leq i\leq q\}$.

We can show that $(H;y_1)\s (H;x_1)$.
For, let $k\geq 0$, and $W\in SW_k(H;y_1)$.
$W$ can be decomposed uniquely to $W=W_1W_2W_3$, 
where $W_1$ and $W_3$ are as long as possible
 and consisting of just $y_1$ and it's  edges in $H$.
Set $f_k(W)=W'_1W_2W'_3$, where $W'_i$ is obtained from $W_i$, by replacing vertex $y_1$ by $x_1$, and edges $y_1z_j$ by $x_1z_j$, for $i=1,3$, and $j=1,2,\cdots,r$. 
Obviously, $f_k(W)\in SW_k(H;x_1)$, and $f_k:SW_k(H;y_1)\to SW_k(H;x_1)$ is injective. Thus $|SW_k(H;y_1)|\leq  |SW_k(H;x_1)|$,
for $k\geq 1$.
 Moreover, $p\geq 2$ implies $|SW_1(H;y_1)|=deg_H(y_1)=r<r+p-1=deg_H(x_1)=|SW_1(H;x_1)|$. Hence $(H;y_1)\s (H;x_1)$.

In a similar method, by changing the end of each semi-edge walk $W\in SW_k(y_i,y_1)$ from $y_1$ to $x_1$, we get that $(H;y_i,y_1)\qs (H;y_i,x_1)$, for $2\leq i\leq q$.

Let $E_{y_1}=\{y_1y_i:2\leq i\leq q\}$, and $E_{x_1}=\{x_1y_i:2\leq i\leq q\}$, and $G=H+E_{x_1}$. 
By lemma \ref{p009},
$SLEE(K_{(p,q)r})=SLEE(H+E_{y_1}) < SLEE(H+E_{x_1})=SLEE(G).$
Note that, since $p\geq 2$, $G$ is a proper subgraph of $K_{(p+q-1,1)r}$. Thus, by lemma \ref{p006},
$SLEE(K_{(p,q)r})<SLEE(G)<K_{(p+q-1,1)r}$
\proofend
By convention, we denote $K_n$ by $K_{(0,1)(n-1)}$, and  $K_{n-1}\cup K_1$ by $K_{(n-1,1)0}$. Now, we can bring in the following results:

\theobegin\label{p013}
Let $G$ be a graph on $n$ vertices, with vertex connectivity $\kappa$, where $0\leq \kappa\leq n-1$. Then $SLEE(G)\leq SLEE(K_{(n-1-\kappa,1)\kappa})$, with equality if an only if $G\cong K_{(n-1-\kappa,1)\kappa}$.
\theoend
\proofbegin
The case $\kappa=n-1$ is trivial, because $K_n$ is the unique graph with vertex (and edge) connectivity $n-1$.

Let $G$ have maximum $SLEE$.
If $\kappa=0$, then the graph is disconnected.  By lemma \ref{p006}, each of it's components must be complete.
 By repeating  use of lemma \ref{p009} (in a similar method used in proof of lemma \ref{p012}), we conclude that $G$ has exactly two components $K_{n-1}$ and $K_1$. Thus $G=K_{n-1}\cup K_1=K_{(n-1,1)0}$.

Now, let $2\leq \kappa \leq n-2$.
Suppose that  $S$ is a subset of $V(G)$, where $G-S$ is disconnected, and
 $|S|=\kappa$.
  By lemma \ref{p006}, $G-S$ is union of two complete  components, say $K_p$ and $K_q$, where $p+q=n-\kappa$.
Again, by lemma \ref{p006}, we have  $G\cong K_{(p,q)\kappa}$.
If $p, q\geq 2$, then lemma \ref{p012} implies that $SLEE(G)<SLEE(K_{(p+q-1,1)\kappa})$, a contradiction. Hence, $q=1$, and $G\cong K_{(n-1-\kappa,1)\kappa}$.
\proofend

Since vertex connectivity of $K_{(n-1-\kappa,1)\kappa}$ is $\kappa$, theorem \ref{p013} ensures that among all graphs with $n$ vertices, $K_{(n-1-\kappa,1)\kappa}$ is the unique graph with maximum $SLEE$.
The following proposition guarantees a similar statement about edge connectivity for $K_{(n-1-\kappa',1)\kappa'}$.
\theobegin\label{p014}
Let $G$ be a graph on $n$ vertices, with edge connectivity $\kappa'$, where $0\leq \kappa'\leq n-1$. Then $SLEE(G)\leq SLEE(K_{(n-1-\kappa',1)\kappa'})$, with equality if an only if $G\cong K_{(n-1-\kappa',1)\kappa'}$.
\theoend
\proofbegin
Suppose that the vertex connectivity of $G$ is $\kappa$.
It is well-known that $\kappa\leq \kappa'$, see \cite{Bondy01}. If $\kappa=\kappa'$, then theorem \ref{p013} implies $SLEE(G)\leq SLEE(K_{(n-1-\kappa',1)\kappa'})$, and equality holds if and only if $G\cong K_{(n-1-\kappa',1)\kappa'}$.
Let $\kappa<\kappa'$. Since $K_{(n-1-\kappa,1)\kappa}$ is a proper subgraph of $K_{(n-1-\kappa',1)\kappa'}$, lemma \ref{p006} and theorem  \ref{p013} yield
$$SLEE(G)\leq SLEE(K_{(n-1-\kappa,1)\kappa})< SLEE(K_{(n-1-\kappa',1)\kappa'})$$
This completes the proof.
\proofend


\begin{thebibliography}{1}
\baselineskip = 0.5cm
\bibitem{Ayyaswamy01}
S.K. Ayyaswamy,
Signless Laplacian Estrada index,
\emph{MATCH Commun. Math. Comput. Chem}
{\bf 66} (2011) 785-794.
\bibitem{Balaban01}
A. T. Balaban, J. Durdevi\'{c}, I. Gutman, 
Comments on $\pi$-electron conjugation in the ﬁve-membered ring of benzo-derivatives of corannulene, 
\emph{Polyc. Arom. Comp.}
{\bf 29}
(2009) 185-205.
\bibitem{Bamdad01}
H. Bamdad, F. Ashraf and Ivan Gutman,
Lower bounds for Estrada index and Laplacian Estrada index,
\emph{Applied Mathematics Letters}
{\bf 23} (2010) 739-742.
\bibitem{Bondy01}
J.A. Bondy and U.S.R Murty,
\emph{Graph Theory with Applications},
American Elsevier, New York, 1976.
\bibitem{Cvetkovic01} 
D. Cvetkovi\'{c}, P. Rowlinson and S. K. Simi\'{c}, 
Signless Laplacians of finite graphs,
\emph{Linear Algebra Appl.}  
{\bf 423}
(2007) 155-171.
\bibitem{Cvetkovic02}
D. Cvetkovi\'{c}, P. Rowlinson, S. K. Simi\'{c}, 
\emph{An Introduction to the Theory of Graph Spectra}, 
Cambridge Univ. Press, Cambridge, 2009.
\bibitem{Das01}
K. C. Das, S. G. Lee,
On the Estrada index conjecture,
\emph{Lin. Algebra Appl.}
{\bf 431} (2009) 1351-1359.
\bibitem{Pena01}
J. A. de la Pe\~{n}a, I. Gutman, J. Rada, 
Estimating the Estrada index, 
\emph{Lin. Algebra Appl.}
{\bf 427} (2007) 70-76.
\bibitem{Deng01}
H. Deng, 
A proof of a conjecture on the Estrada index, 
\emph{MATCH Commun. Math. Comput. Chem.}
{\bf 62} (2009) 599–606.
\bibitem{Dias01}
J. R. Dias,
Structure/formula informatics of isomeric sets of ﬂuorantheno-
id/ﬂuorenoid and indacenoid hydrocarbons, 
\emph{J. Math. Chem.}
{\bf 48} 
(2010) 313-329.
\bibitem{Du01}
Z. Du and B. Zhou,
The Estrada index of trees,
\emph{Linear Algebra Appl.} 
{\bf 435 }
(2011) 2462-2467.
\bibitem{Du02}
Z.Du, B. Zhou and R. Xing,
On maximum Estrada indices of graphs with given parameters,
\emph{Linear Algebra and its Applications}
{\bf 436} (2012) 3767-3772.
\bibitem{Estrada01}
E. Estrada,
Characterization of 3D molecular structure,
\emph{ Chem. Phys. Lett.} 
{\bf 319 }
(2000) 713-718. 
\bibitem{Estrada02}
E. Estrada, 
Characterization of the folding degree of proteins, 
\emph{Bioinformatics} 
{\bf 18 }
(2002) 697-704.
\bibitem{Estrada03}
E. Estrada, 
Characterization of the amino acid contribution to the folding degree of proteins, 
\emph{Proteins} 
{\bf 54 }
(2004) 727-737.
\bibitem{Estrada04}
E. Estrada, N. Hatano, 
Statistical-mechanical approach to subgraph centrality in complex networks, 
\emph{Chem. Phys. Lett.} 
{\bf 439 }
(2007) 247-251.
\bibitem{Estrada05}
E. Estrada, J. A. Rodr\'{i}guez-Vel\'{a}zquez, 
Subgraph centrality in complex networks,
\emph{Phys. Rev. E} 
{\bf 71 }
(2005) 056103-1-056103-9.
\bibitem{Estrada06}
E. Estrada, J. A. Rodr\'{i}guez-Vel\'{a}zquez, 
Spectral measures of bipartivity in complex
networks, 
\emph{Phys. Rev. E} 
{\bf 72 }
(2005) 046105-1-046105-6.
\bibitem{Estrada07}
E. Estrada, J. A. Rodr\'{i}guez-Vel\'{a}zquez, M. Randi\'{c}, 
Atomic branching in molecules,
\emph{Int. J. Quantum Chem. }
{\bf 106}
(2006) 823-832.
\bibitem{FathTabar01}
 G.H. Fath-Tabar, A.R. Ashrafi, and I. Gutman, 
Note on Estrada and $L$-Estrada indices of graphs,
\emph{Bull. Cl. Sci. Math. Nat. Sci. Math.} 
{\bf 34 }
(2009) 1-16.
\bibitem{Grone01}
R. Grone, R. Merris,
The Laplacian spectrum of a graph II, 
\emph{SIAM J. Discr. Math.}
{\bf 7} 
(1994) 221-229.
\bibitem{Grone02} 
R. Grone, R. Merris, V. S. Sunder, 
The Laplacian spectrum of a graph, 
\emph{SIAM J. Matrix Anal. Appl.}
{\bf 11} 
(1990) 218-238.
\bibitem{Gutman04} 
I. Gutman, S. Radenkovi\'{c}, 
A lower bound for the Estrada index of bipartite molecular graphs,
\emph{Kragujevac J. Sci. }
{\bf 29} (2007) 67-72.
\bibitem{Gutman05} 
I. Gutman, E. Estrada, J. A. Rodr\'{i}guez-Vel\'{a}zquez, 
On a graph-spectrum-based structure descriptor, 
\emph{Croat. Chem. Acta }
{\bf 80} (2007) 151-154.


\bibitem{Gutman01}
I. Gutman, J. Durdevi\'{c},
Fluoranthene and its congeners - A graph theoretical study, 
\emph{MATCH Commun. Math. Comput. Chem.}
{\bf 60} (2008) 659-670.
\bibitem{Gutman02}
I. Gutman, J. Durdevi\'{c}, 
On $\pi$-electron conjugation in the ﬁve-membered ring of ﬂuoranthene-type benzenoid hydrocarbons, 
\emph{J. Serb. Chem. Soc. }
{\bf 74} (2009) 765-771.
\bibitem{Gutman03}
I. Gutman, J. Durdevi\'{c},
Cycles in dicyclopenta-derivatives of benzenoid hydro-carbons, 
\emph{MATCH Commun. Math. Comput. Chem.}
{\bf 65} (2011) 785-798.
\bibitem{Gutman06}
I. Gutman,
Lower bounds for Estrada index,
\emph{Publ. Inst. Math.(Beograd)}
{\bf 83} (2008) 1-7.
\bibitem{Khosravanirad01}
A. Khosravanirad,
A lower bound for Laplacian Estrada index of a graph,
\emph{MATCH Commun. Math. Comput. Chem}
{\bf 70} (2013) 175-180.
\bibitem{Li01}
J. Li, W.C. Shiu and W.H. Chan,
Note on the Laplacian Estrada index of a graph,
\emph{MATCH Commun. Math. Comput. Chem}
{\bf 66} (2011) 777-784.
\bibitem{Merris01}
R. Merris, 
Laplacian matrices of graphs: a survey, 
\emph{Linear Algebra Appl.}
{\bf 197-198} (1994) 143-176.
\bibitem{Zhang01}
J. Zhang, B. Zhou and J. Li,
On Estrada index of trees,
\emph{Linear Algebra and its Applications},
{\bf 434} (2011) 215-223.
\bibitem{Zhao01}
H. Zhao and Y. Jia,
On the Estrada index of bipartite graph,
\emph{Communications in Mathematical and in Computer Chemistry/MATCH }
{\bf 61} (2009) 495-501.
\bibitem{Zhou01}
B. Zhou and I. Gutman,
More on the Laplacian Estrada Index,
\emph{Applicable Analysis and Discrete Mathematics}
{\bf 3} (2009) 371-378.
\end{thebibliography}
\end{document}